\documentclass[review]{elsarticle}

\usepackage{lineno,hyperref, amsmath,amsthm,amsfonts, amssymb, euscript, enumitem}
\modulolinenumbers[5]

\journal{Journal of \LaTeX\ Templates}

\theoremstyle{definition}

\begin{document}

\begin{frontmatter}

\title{Metrically Ramsey ultrafilters}

\author{Igor Protasov and Ksenia Protasova }

\begin{abstract}
Given a metric space $(X,d)$, we say that a mapping $\chi: [X]^{2}\longrightarrow\{0.1\}$ is an isometric coloring if $d(x,y)=d(z,t)$ implies $\chi(\{x,y\})=\chi(\{z,t\})$.  A free ultrafilter $\mathcal{U}$
on an infinite metric space $(X,d)$ is called  metrically Ramsey  if, for every isometric coloring $\chi$ of $[X]^{2}$, there is a member $U\in\mathcal{U}$ such that  the set  $[U]^{2}$ is $\chi$-monochrome. We prove that each infinite ultrametric space $(X,d)$ has a countable subset $Y$ such that each free ultrafilter $\mathcal{U}$ on $X$ satisfying $Y\in\mathcal{U}$ is metrically Ramsey. On the other hand, it is an open question whether every metrically Ramsey ultrafilter on the natural numbers $\mathbb{N}$ with the metric $|x-y|$  is a Ramsey ultrafilter. We prove that every metrically Ramsey ultrafilter $\mathcal{U}$ on $\mathbb{N}$ has a member with no arithmetic progression of length 2, and if $\mathcal{U}$ has a thin member then there is a mapping $f:\mathbb{N}\longrightarrow\omega $ such that $f(\mathcal{U})$  is a Ramsey ultrafilter.

{\bf Classification:}  03E05, 54E35


\end{abstract}

\end{frontmatter}


\large

Let $X$ be an infinite set and let $\mathfrak{F}$ be some family of $\{0,1\}$-colorings of the set $[X]^{2}$ of all two-element subsets of $X$. 
We say that a free   ultrafilter $\mathcal{U}$ on $X$ is {\it  Ramsey with respect to } $\mathfrak{F}$ if, for any coloring $\chi\in\mathfrak{F}$, there exists $U\in\mathcal{U}$ such that $[U]^{2}$ is $\chi$-monochrome.
In the case in which $\mathfrak{F}$ is the family of all $\{0,1\}$-colorings of  $[X]^{2}$, we get the classical definition of {\it  Ramsey ultrafilters}. 
It is well-known that $\mathcal{U}$ is a Ramsey ultrafilter if and only if $\mathcal{U}$ is {\it selective}, i.e. for every partition $\mathcal{P}$ of $X$
either  $P\in\mathcal{U}$ for some $P\in\mathcal{P}$ or there exists $U\in\mathcal{U}$ such that $|U\bigcap P|\leq 1$  for each $P\in\mathcal{P}$.

Given a metric space $(X,d)$, we say that a mapping $\chi: [X]^{2}\longrightarrow\{0,1\}$ is an {\it isometric coloring} if $d(x,y)=d(z,t)$ implies $\chi(\{x,y\})=\chi(\{z,t\})$. We note that every isometric coloring $\chi$ is uniquely defined by some mapping $f: d(X,X)\setminus \{0\} \longrightarrow \{0,1\}$. Indeed, we take an arbitrary $r\in d(X,X)\setminus \{0\}$, choose $\{x,y\}\in [X]^{2}$ such that $d(x,y)=r$ and put $f(r)= \chi(\{x,y\})$.
On the other hand, for
$f: d(X,X)\setminus \{0\} \longrightarrow \{0,1\}$, we define $\chi$ by $\chi(\{x,y\})=f(d(x,y))$.

We say that a free  ultrafilter on an infinite metric space $(X, d)$ is {\it  metrically Ramsey} if $\mathcal{U}$ is
Ramsey with respect to all isometric colorings of $[X]^{2}$.

Let $G$ be a group and let $X$ be a $G$-space with the action $(G,X)\longrightarrow X$, $(g,x)\longmapsto gx$. A coloring  $\chi: [X]^{2}\longrightarrow\{0.1\}$ is called {\it $G$-invariant}  if  $\chi(\{x,y\})=\chi(\{gx,gy\})$ for all $\{x,y\}\in [X]^{2}$ and $g\in G$. A free ultrafilter $\mathcal{U}$ of $X$ is called {\it $G$-Ramsey} if $\mathcal{U}$ is Ramsey with respect to the family of all $G$-invariant colorings of $[X]^{2}$.

We consider the special case: $X$  is a metric space and $G$  is a group of isometries of $X$. Clearly, every isometric coloring of  $[X]^{2}$ is $G$-invariant. If $G$ is metrically 2-transitive (if $d(x,y)=d(z,t)$ then there is $g\in G$ such that $g\{x,y\}=\{z,t\}$ ) then every $G$-invariant coloring of $[X]^{2}$ is an isometric coloring.

We take the group  $\mathbb{Z}$ of integers, put  $X=\mathbb{Z}$ and consider the action  $\mathbb{Z}$  on $X$  by $(g,x)= g+x$. {\it Is every  $\mathbb{Z}$-Ramsey ultrafilter selective?}
This question appeared in \cite{b5} and, to our knowledge, remains open. We endow  $\mathbb{Z}$ with the metric $d(x,y)=|x-y|$. By above paragraph an ultrafilter $\mathcal{U}$ on $\mathbb{Z}$ is $\mathbb{Z}$-Ramsey if and only if  $\mathcal{U}$ is metrically Ramsey.
{\it Is every metrically Ramsey ultrafilter on $\mathbb{Z}$-selective?}
This is an equivalent form of the above question. The case of
$\mathbb{Z}$ evidently equivalent to the  case of $\mathbb{N}$.

Surprisingly or not, the case of ultrametric spaces is cardinally different and much more easy  to explore. We recall that a metric $d$ is an {\it ultrametric} if $d(x,y)\leq \max\{d(x,z), d(z,y)\}$ for all  $x,y,z\in X$. We prove that every infinite ultrametric space $X$ has a countable subset $Y$ such that any ultrafilter $\mathcal{U}$ on $X$  satisfying  $Y\in U$ is metrically Ramsey.

\section{ \large Equidistance subsets}

We say that a subset $Y$ of a metric space $(X,d)$ is an {\it
equidistance subset} if there is $r\in\mathbb{R}^{+}$ such that $d(x,y)=r$ for all distinct $x,y\in Y$. If $Y$ is an equidistance subset of $(X,d)$ then every free ultrafilter $\mathcal{U}$ on $X$ such that $Y\in \mathcal{U}$  is  metrically Ramsey.
\vspace{3 mm}

{\bf Propozition 1.1.} {\it Every infinite metric space with finite scale $d(X,X)$, $d(X,X)=\{d(x,y): x,y\in X\}$ has a countable equidistance subset.}
\vspace{1 mm}

\begin{proof} We define a coloring $\chi : [X]^{2}\longrightarrow d(X,X)$ by $\chi (\{x,y\})= d(x,y)$ and apply the classical Ramsey theorem [2, p.16].
\end{proof}

For an ultrametric space $(X,d)$ and $r\in d(X,X)$, we use the equivalence
 $  \   \sim _{r} \  $
 defined by $x \sim _{r}  y$ if and only if $d(x,y)\leq r$. Then $X$ is partitioned into classes of $r$-equivalence
$X=\bigcup_{\alpha<\lambda} X_{\alpha}$. If $x,y\in X_{\alpha}$ then $d(x,y)\leq r$. If $x\in X_{\alpha}$, $y\in X_{\beta}$ and $\alpha \neq\beta$ then $d(x,y)>r$.

\vspace{3 mm}

{\bf Propozition 1.2.} {\it Let $(X,\alpha)$  be an infinite ultrametric space with finite scale $d(X,X)$. If $|X|$ is regular then $X$ has an equidistance subset $Y$ of cardinality $|Y|=|X|$. If $|X|$ is singular then, for every cardinal $\kappa<|X|$ there is an equidistance subset of cardinality $\kappa$.}
\vspace{1 mm}

\begin{proof}
Let $d(X,X)= \{0, r_{1}, \ldots , r_{n}\}$,  $ 0<r_{1}< \ldots < r_{n} $. We proceed on induction by $n$. For $n=1$, the statement is evident: $Y=X$.

To make the inductive step from $n$  to  $n+1$, we partition $X$ into classes of $r_{n}$-equivalence $X=\bigcup_{\alpha<\lambda} X_{\alpha}$. If $\lambda=|X|$  then  we pick one element $y_{\alpha}\in  X_{\alpha}$ and put $Y=\{y_{\alpha}: \alpha< |X|\}$, so $d(x,y)= r_{n+1}$ for all distinct $x,y\in Y$. Assume that $\lambda< |X|$. If $|X|$ is regular, we take $\alpha$ so that $|X_{\alpha}|=X$ and apply the inductive assumption to $X_{\alpha}$. If $|X|$ is singular then we take $X_{\alpha}$ such that $|X_{\alpha}|>\kappa$ and use the inductive assumption.
\end{proof}

\vspace{1 mm}

{\bf Remark 1.1.} If in  Proposition 1.2. $|X|$ is singular, we cannot state that there is an equidistance subset $Y$ of cardinality $X$. We take an arbitrary singular cardinal $\kappa$, put $X=\kappa$ and partition $X=\bigcup_{\alpha<\lambda} X_{\alpha}$ so that $\lambda<\kappa$ and
$|X_{\alpha}|< \kappa$ for each $\alpha<\kappa$.  We define an ultrametric $d$ on $X$ by $d(x,x)=0$, $d(x,y)=1$ if $x,y\in X_{\alpha}$, $x\neq y \  $ and $ \ d(x,y)=2$ if $x\in X_{\alpha}, \   \  y\in X_{\beta}$, $ \  \alpha\neq\beta$.  if $Y$
is an equidistance subset of $(X,\alpha)$ then either $Y\subseteq X_{\alpha}$ for some $\alpha<\lambda$, or $|Y\bigcap X_{\alpha}|\leq 1$ for each $\alpha<\lambda$. Hence, $|Y|<\kappa$.
\vspace{3 mm}

{\bf Proposition 1.3.} {\it For every infinite cardinal $\kappa$, there exists a metric space $(X,d)$ such that $|X|=2^{\kappa}$, $d(X,X)=\{0,1,2\}$  and every equidistance subset $Y$ of $(X,d)$ is of cardinality $|Y|\leq\kappa$.}

\begin{proof} We put $X=2$  and apply [4, Theorem 6.2] to define a coloring $\chi: [X]^{2}\longrightarrow\{1,2\}$ with no monochrome $[Z]^{2}$ for $|Z|>\kappa$. Then we define a metric $d$ on $X$
 by $d(x,x)=0$ and $d(x,y)=\chi(x,y)$ for all distinct $x,y\in X$.
 \end{proof}
\vspace{1 mm}

{\bf Propozition 1.4.} {\it Let $(X,d)$ be a metric space with infinite scale  $d(X,X)$, $|d(X,X)|=\kappa$.
If $|X|\geq (2^{\kappa})^{+}$  then
there is an equidistance subset $Y$  of $(X, d)$ such that $|Y|=\kappa^{+}$.}

\begin{proof} We define a coloring
$\chi: [X]^{2}\longrightarrow  d(X,X)$ by $\chi(\{x,y\})=d(x,y)$ and apply the Erd$\acute{o}$s-Rado theorem [4, Theorem 6.4].
 \end{proof}
\vspace{1 mm}

{\bf Proposition 1.5.} {\it For every infinite
 metric space $(X,d)$,  there exists an injective sequence $(x_{n})_{n\in\omega}$ in $X$ such that one of the following conditions is satisfied:

 $(i)$ the sequence  $(d(x_{0}, x_{n})_{n\in\omega}$ is increasing;

 $(ii)$ the sequence  $(d(x_{0}, x_{n})_{n\in\omega}$ is decreasing;

 $(iii)$ for every $n\in\omega \  $ and all $  \  i,j $,  $ \ i>n$, $  \  j>n$, we have $  \ d(x_{n}, x_{i})=d(x_{n}, x_{j})$.}

\begin{proof} We assume that there exists $x_{0}\in X$ such that the set $d(x_{0}, X)$ is infinite, $d(x_{0}, X)= \{d(x_{0}, x) :x\in X\}$. We choose a countable subset $Y$  of $X$  such that $x_{0}\notin Y$ and $d(x_{0}, y)\neq d(x_{0}, z)$ for all distinct $y,z\in Y$.  The set $d(x_{0}, Y)$  contains either increasing or decreasing sequence $(r_{n+1})_{n\in\omega}$. For each $n\in\omega$, we choose $x_{n+1}$ such that $d(x_{0}, x_{n+1})=r_{n}$. Then the sequence $(d(x_{0}, x_{n}))_{n\in\omega}$ satisfies either $(i)$ or $(ii)$.

In the alternative case, the set $d(x, X)$  is finite for each $x\in X$. We fix $x_{0}\in X$ and choose a countable subset $X_{1}$ such that $|d(x_{0}, X_{1})|=1$. We pick $x_{1}\in X_{1}$ and choose a countable subset $X_{2}\subseteq X_{1}$ such that$|d(x_{1}, X_{2})|=1$ and so on. After $\omega$ steps, we get the sequence $(x_{n})_{n\in\omega}$ satisfying  $(iii)$.
 \end{proof}
 \vspace{1 mm}

 {\bf Proposition 1.6.} {\it Let $(X,d)$ be an infinite metric space and let $\{x_{n} : n\in\omega \}$ be a family of  non-empty pairwise disjoint subsets of $X$ such that $d(X_{i}, X_{j}) \bigcap d(X_{n}, X_{n})=\emptyset$ for all $n$  and distinct $i, j$. Let $\mathcal{U}$ be a metrically
Ramsey ultrafilter on $X$ such that $\bigcup_{n<\omega} X_{n}\in\mathcal{U}$ and $X_{n}\notin\mathcal{U}$
for each $n<\omega$.  Then the following statements hold:
\vspace{3 mm}

 $(i)$  there exists $U\in\mathcal{U}$ such that $|U\bigcap X_{n}|\leq 1$ for each $n<\omega$;
 \vspace{3 mm}

 $(ii)$ if $d(X_{i}, X_{j}) \bigcap d(X_{k}, X_{l})=\emptyset$ for all distinct $\{i,j\},  \{k,l\}\in[\omega]^{2}$ then there is a mapping $\varphi : X\longrightarrow\omega$ such that the ultrafilter $\varphi(\mathcal{U})$ is selective;
 \vspace{3 mm}

$(iii)$ if for each $n<\omega$ there exists $m<\omega$ such that $m>n$, $|X_{m}|>n$, then there exists an ultrafilter $\mathcal{V}$ on $X$ such that $\bigcup_{n<\omega}\in\mathcal{V}$ and $\mathcal{V}$ is not metrically Ramsey. }
\vspace{2 mm}

\begin{proof} $(i)$ By the assumption, the sets $A=\bigcup_{n<\omega}d(X_{n}, X_{n})$ and $B=\bigcup_{i\neq j}\ d(X_{i}, X_{j}) $ are disjoint. We take an arbitrary mapping $f: \mathbb{R}^{+}\longrightarrow\{0,1\}$ such that $f|_{A} \equiv 0$, $f|_{B} \equiv 1$, and consider the isometric coloring $\chi$ of $[X]^{2}$ defined by $f$. Since $\mathcal{U}$ is metrically Ramsey, there is $U\in \mathcal{U}$ such that $[U]^{2}$ is $\chi$-monochrome. Clearly, $|U\bigcap X_{n}|\leq 1$  for each $n<\omega$.\vspace{3 mm}

$(ii)$ We define $\varphi$ by the rule: if $x\in X_{i}$  then $\varphi(x)=i$, if $x\in X\setminus \bigcup_{n<\omega}X_{n}$ then $\varphi(x)=0$. We take an arbitrary coloring $\chi^{\prime}: [\omega]^{2}\longrightarrow\{0,1\}$ and define a coloring $\chi: [\bigcup_{n<\omega} X_{n}]^{2}\longrightarrow\{0,1\}$ as follows. If $x\in X_{i}$,  $y\in X_{j}$, $i\neq j$ then $\chi (\{x,y\})= \chi^{\prime} (\{x,y\})$. If $x,y \in X_{n}$ then $\chi (\{x,y\})=0$. By the assumption, the coloring $\chi$ is isometric. We choose $U\in\mathcal{U}$ such that $U\subseteq \bigcup_{n<\omega} X_{n}$ and $[U]^{2}$ is $\chi$-monochrome and $|U\bigcap X_{n}|\leq 1$ for each $n<\omega$. Then $\varphi(U)\in \varphi(\mathcal{U})$ and $[\varphi(U)]$ is $\chi^{\prime}$-monochrome, so $\varphi(\mathcal{U})$ is a Ramsey ultrafilter.\vspace{3 mm}

$(iii)$ We consider the family of all filters $\mathfrak{F}$ on $X$  such that $\bigcup_{n<\omega}X_{n}\in\mathfrak{F}$  and, for every $n\in\omega$ and $F\in\mathfrak{F}$, there exists $m\in\omega$ such that $|F\bigcap X_{m}|> n$. By the Zorn Lemma, this family has maximal by inclusion element $\mathcal{V}$. It is easy to verify that $\mathcal{V}$ is ultrafilter. By  $(i)$,  $\mathcal{V}$ is not metrically Ramsey.
 \end{proof}

\section{\large The ultrametric case}

 {\bf Proposition 2.1.} {\it For every infinite ultrametric space $(X, d)$, there exists a countable subset $Y$  of $X$  such that every free ultrafilter $\mathcal{U}$ on $X$ satisfying $Y\in\mathcal{U}$ is metrically Ramsey.}
  \vspace{2 mm}

  \begin{proof} We choose the sequence $(x_{n})_{n\in\omega}$ given by Proposition 1.5, put $Y=\{x_{n}: n\in\omega\}$, fix an arbitrary mapping $f: \mathbb{R}^{+}\longrightarrow\{0,1\}$ and take an arbitrary free ultrafilter $\mathcal{U}$ satisfying $Y\in \mathcal{U}$.

  We assume that either $(i)$ or $(ii)$   of Proposition 1.5 hold for $(x_{n})_{n\in\omega}$. We define a mapping $h: Y\longrightarrow \mathbb{R} ^{+}$ by $h(x_{n})=d(x_{0}, x_{n})$  and choose $k\in\{0,1\}$ such that $(fh)^{-1}(k)\in\mathcal{U}$. Since $d$ is an ultrametric, in the case
  $(i)$ we have $d(x_{i}, x_{n})= d(x_{0}, x_{n})$ for all $i< n $, and in the case $(ii)$ we have
 $d(x_{i}, x_{n})= d(x_{0}, x_{i})$ for all $i< n $. In  both cases, if $\{x_{i}, x_{n}\}\in [(fh)^{-1}(k)]^{2}$ then $f(d(x_{i}, x_{n}))= k$.

 If $(x_{n})_{n\in\omega}$ satisfies  $(iii)$ of Proposition 1.5 then we define a mapping  $h: Y\longrightarrow \mathbb{R} ^{+}$ by $h(x_{n})=d(x_{n}, x_{i})$, $i>n$ and repeat above arguments.
  \end{proof}
 \vspace{2 mm}

 {\bf Proposition 2.2.} {\it For a free ultrafilter $\mathcal{U}$ on an infinite set $X$, the following statements are equivalent:
 \vspace{3 mm}

 $(i)$  $\mathcal{U}$ is selective;
\vspace{3 mm}

$(ii)$  $\mathcal{U}$ is metrically Ramsey for each ultrametric $d$ on $X$  such that  $d(X,  X)= \{0,1,2\}$.}
 \vspace{3 mm}

\begin{proof} The implication  $(i) \Longrightarrow$  $(ii)$ is evident. We assume that $\mathcal{U}$ is not selective and choose a partition  $\mathcal{P}$ of $X$ such that $P\notin\mathcal{U}$ for each $P\in\mathcal{P}$, and for every $U\in\mathcal{U}$, there is $P\in\mathcal{P}$  such that $|P\bigcap U|>1$. We define an ultrametric $d$ on $X$ by $d(x,x)=0$,  $d(x,y)=1$  if $x\neq y$, $x,y\in P$ for some $P\in\mathcal{P}$, and $d(x,y)=2$ if $x,y$ belong to different cells of the partition $\mathcal{P}$.   We define a coloring $\chi: [X]^{2}\longrightarrow\{1,2\}$ by $\chi(\{x,y\})$ $ = d(x,y)$. Then the set $[U]^{2}$ is not $\chi$-monochrome for each $U\in\mathcal{U}$  so  $\mathcal{U}$ is not metrically Ramsey and $(ii) \Longrightarrow$  $(i)$.
\end{proof}
 \vspace{1 mm}

 {\bf Proposition 2.3.} {\it Let $(X, d)$  be an infinite ultrametric space with finite scale  $d(X,  X)= \{0,r_{1},\ldots , r_{n}\}$, $0< r_{1}< \ldots < r_{n} $. Then the following statements are equivalent:

 \vspace{3 mm}

 $(i)$  every free ultrafilter on $(X,d)$ is metrically Ramsey;

 \vspace{3 mm}

 $(ii)$  for every $i\in \{1, \ldots , n\}$, the partition $\mathcal{P}_{i}$ of $X$ into classes of $r_{i}$-equivalence has only finite number of infinite classes and there is $m\in\omega$ such that $|C|<m$ for each finite class $C$ from $\mathcal{P}_{i}$.}

\begin{proof}  $(i) \Longrightarrow$  $(ii)$. If $\mathcal{P}_{i}$  has infinitely many infinite classes or the set $ \{|C|: C $ is a finite class from $\mathcal{P}_{i}\}$ is infinite we apply Proposition $1.6 (iii)$ to get a free ultrafilter $\mathcal{V}$ on $X$ which is not metrically Ramsey.

$(ii) \Longrightarrow$  $(i)$. We proceed on induction by $n$. For $n=1$, the statement is evident.

 \vspace{2 mm}

 To make the inductive step from $n$  to $n+1$, take an arbitrary free ultrafilter $\mathcal{U}$  on $X$ and  we consider the partition $\mathcal{P}_{n+1}$. 
 Let  $X_{1},\ldots , X_{m} \ $ be the set of all infinite classes from  $ \ \mathcal{P}_{n+1}$. 
 If $ \ X_{1}\bigcup$ $\ldots \bigcup X_{m}\in \mathcal{U}$  then we take $X_{i}\in\mathcal{U}$ and apply the inductive assumption.
 If  $X\setminus (X_{1}\bigcup  \ldots  \bigcup X_{m})\in \mathcal{U}$  then we choose $U\in\mathcal{U}$ such that $U\subseteq X\setminus(X_{1}\bigcup , \ldots , \bigcup X_{m})$  and  $|U\bigcap C|\leq 1$ for each finite class
 $C\in\mathcal{P}_{n+1}$. Then $U$ is an equidistance set so $\mathcal{U}$ is metrically Ramsey.
\end{proof}

\section{\large The case of $\mathbb{N}$}

{\bf Proposition 3.1.} {\it Let $\mathcal{U}$  be a metrically Ramsey
ultrafilter  on $\mathbb{N}$ and let $f: \mathbb{N}\longrightarrow  N$ be a mapping such that $f(x)>x$ for each
$x\in\mathbb{N}$. Then there exists a member $U\in\mathcal{U}$ having no subsets of the form $\{a, a+x, a+ f(x)\}$.  In particular $(f(x)=2x)$, some member of $\mathcal{U}$  has no arithmetic progressions of length 2. }\vspace{3 mm}

\begin{proof} We consider a directed graph $\Gamma_{f}$ with the set of vertices $\mathbb{N}$ and the set of edges
 $\{ (x,f(x)) : \ x\in\mathbb{N}\}$. Since $f(x)>x$, $\Gamma_{f}$ is the disjoint union of directed trees $T$ such that each vertex of $T$ has at most one input edges. Using  this observation, it is easy to partition $\mathbb{N}= A_{1}\bigcup A_{2}$ so that $f(A_{1})\subseteq A_{2}$,  $ \ f(A_{2})\subseteq A_{2}$.

 The partition $N=A_{1}\bigcup A_{2}$ defines an isometric coloring $\chi : [\mathbb{N}]^{2}\longrightarrow\{1,2\}$ by $\chi(\{x,y\})=i$ if and only if $d(x,y)\in A_{i}$. We take a subset $U\in\mathcal{U}$ such that the set $[U]^{2}$ is $\chi$-monochrome and assume that $\{a, a+x, a+ f(x)\}\subset U$ for some $a,x\in \mathbb{N}$. We note that $d(a,a+x)=x$,
 $ \ d(a, a+f(x))= f(x)$, but $x$  and $f(x)$ belong to different subsets $A_{1}, A_{2}$, so
  $\chi( \{a,a+x\})\neq \chi( \{a,a+f(x)\})$ and  we get  a contradiction with the choice of $U$.
 \end{proof}\vspace{2 mm}

Let $\mathcal{U}$ be metrically Ramsey ultrafilter  on $\mathbb{N}$. Assume that there is $U\in \mathcal{U}$  such that $d(x,y)\neq d(z,t)$ for all distinct $\{x,y\}, \{z, t\}\in [U]^{2}$. Then every $\{0,1\}$-coloring of $[U]^{2}$ can be extended to some isometric coloring of $[\mathbb{N}]^{2}$.  Hence, $\mathcal{U}$  is a  Ramsey ultrafilter.

We say that a subset $T=\{t_{n}: t_{n}< t_{n+1},  n<\omega\}$  of $\mathbb{N}$ is {\it thin} if $(t_{n+1}- t_{n})\longrightarrow\infty \ $ as $ \ n\longrightarrow\infty$. \vspace{3 mm}


{\bf Proposition 3.2.} {\it If a metrically Ramsey ultrafilter $\mathcal{U}$   on $\mathbb{N}$ has a thin subset $T\in\mathcal{U}$ then there exists a mapping $\varphi: \mathbb{N}\longrightarrow\omega$ such that the ultrafilter $\varphi(\mathcal{U})$ is selective and $\varphi$ is finite-to-one on some member $U\in\mathcal{U}$.}\vspace{2 mm}

\begin{proof}
Let $T=\{t_{n}: t_{n}< t_{n+1}, n\in\omega\}$.  Assume that we have chosen two sequences $(a_{n})_{n\in\omega},  \ (b_{n})_{n\in\omega} \ $ in $T$  such  that \\
$(1) \ $  $ \ a_{n}< b_{n}< a_{n+1}< b_{n+1} \ $ for each $  \  n\in\omega$; \\
$(2) \ $   $ \ d([a_{n}, b_{n}) \ \bigcap \ T, \ [a_{n}, b_{n}) \ \bigcap \  T) \ \bigcap  \ d([a_{i}, b_{i}) \ \bigcap \ T, \  [a_{j}, b_{j}) \  \bigcap \  T) \ = \ \emptyset$
for all $n$ and distinct $i, j$ ; \\
$(3)$   $ \ d([a_{i}, b_{i})  \bigcap T,  [a_{j}, b_{j})  \bigcap   T)  \bigcap   d([a_{k}, b_{k})  \bigcap  T,   [a_{l}, b_{l})   \bigcap  T)  =  \emptyset$
for all distinct $\{i, j\}$,  $ \ \{k, l\}\in [\omega]^{2}$ ; \\
$(4)$   $ \ d([b_{n}, a_{n+1}) \bigcap  T,  [b_{n}, a_{n+1})  \bigcap   T)  \bigcap   d([b_{i}, a_{i+1})  \bigcap  T,   [b_{j}, a_{j+l})   \bigcap  T)  = \emptyset$ for all $n$ and distinct $i, j$ ; \\
$(5)$   $ \ d([b_{i}, a_{i+1})  \bigcap  T,  [b_{j}, a_{j+1})  \bigcap   T) \bigcap  d([b_{k}, a_{k+1}) \ \bigcap T,   [b_{l}, a_{l+l})   \bigcap  T)  =  \emptyset$
for all distinct $\{i, j\}$,  $ \ \{k, l\}\in [\omega]^{2}$ ;
\vspace{2 mm}

We put $A= \ \bigcup_{n\in\omega}
([a_{n},b_{n})  \bigcap  T)$,  $ \ B=\bigcup_{n\in\omega} ([b_{n}, a_{n+1})  \bigcap   T)$ and note that $A$ and $B$  with corresponding partitions satisfy Proposition 1.6 $(ii)$. Since either $A\in \mathcal{U}$ or $B\in \mathcal{U}$,
Proposition 1.6 $(ii)$ gives the mapping $\varphi: \mathbb{N}\longrightarrow\omega$ such that $\varphi(\mathcal{U})$ is selective.  By  the construction of $\varphi$, $ \ \varphi$ is finite-to-one on $A$ or $B$ respectively.

It remains to construct $(a_{n})_{n\in\varphi}$ and $(b_{n})_{n\in\varphi}$. We put $a_{0}= t_{0}$, $ \ b_{0}= t_{1}$ and assume that we have chosen $a_{0}, b_{0}, \ldots , a_{n}, b_{n}$. Since $T$  is thin, we can choose $a_{n+1}\in T$ so that $a_{n+1}> 2b_{n}$ and $|t-t^{\prime}|> 2 a_{n}$ for all distinct $t,t^{\prime}\in T\setminus [1, a_{n+1})$. Then we choose $b_{n+1}\in T$ so that $b_{n+1}> 2  a_{n+1} $ and $|t-t^{\prime}|> 2 b_{n}$ for all distinct $t,t^{\prime} \in T\setminus [1, b_{n+1})$.  After $\omega$ steps, we get the desired $(a_{n})_{n\in\varphi} (b_{n})_{n\in\varphi}$.
\end{proof}

\section{\large Comments and open questions}

1. In connection with Proposition 3.1, we mention [5, Corollary 2]: every metrically Ramsey ultrafilter on $X$  has a member $U$ with no subsets of the form $\{x, y, x+y\}$, $x\neq y$.
\vspace{3 mm}

In connection with Proposition 3.2, we ask

\vspace{5 mm}
{\bf Question 4.1.} {\it Let $\mathcal{U}$ be a metrically Ramsey ultrafilter on $\mathbb{N}$. Does there exist a thin subset $U\in \mathcal{U}$?}

\vspace{5 mm}
{\bf Question 4.2.} {\it Assume that a metrically Ramsey ultrafilter  $\mathcal{U}$  on $\mathbb{N}$ has a thin member. Is  $\mathcal{U}$ selective?}
\vspace{5 mm}

 2. Let  $G$  be an Abelian  group. A coloring $\chi: [G]^{2}\longrightarrow\{0,1\}$ is called a PS-{\it coloring} if, for $\{x,y\}, \{z,t\}\in [G]^{2}$,  $x+y=z+t$ implies $\chi(\{x,y\})=\chi(\{z,t\})$. A free ultrafilter
 $\mathcal{U}$ on $G$ is called a PS-{\it ultrafilter } if
 $\mathcal{U}$ is Ramsey  with respect to all
 PS-{\it colorings} of $[G]^{2}$. The PS-ultrafilters were introduced and studied in \cite{b5}, for exposition of \cite{b6} see [1, Chapter 10].

 If $G$  has a finite set $B(G)=\{g\in G:2g=0\}$ of elements of order 2 then every PS-ultrafilter on $G$ is selective. If $B(G)$ is infinite then, under Martin's Axiom, there is a non-selective  PS-ultrafilter on $G$. If there exists  PS-ultrafilter on some countable group $G$ then there is a $P$-point in $\omega^{\ast}$.

 Now we consider the countable Boolean group $\mathbb{B}$,  $  \  B(\mathbb{B})=\mathbb{B}$. We note that a coloring
 $\chi: [\mathbb{B}]^{2}\longrightarrow\{0,1\}$ is a PS-coloring if and only  if $\chi$ is $\mathbb{B}$-invariant.
 Thus, a free ultrafilter  $\mathcal{U}$ on $\mathbb{B}$ is a PS-ultrafilter if and only if $\mathcal{U}$ is  $\mathbb{B}$-Ramsey. By above paragraph, in the models of ZFC with no P-points in $\omega^{\ast}$, there are no
 $\mathbb{B}$-Ramsey ultrafilters. 
 However, every strongly summable ultrafilter on $\mathbb{B}$ is a $PS$-ultrafilter. For  strongly summable ultrafilters on Abelian groups see \cite{b3}.

On the other hand, $\mathbb{B}$ is the direct $ sum\oplus_{n<\omega} \ \{0,1\}_{n}$
of $\omega$ copies of $\mathbb{Z}_{2}= \{0,1\}$,
and has the natural structure of ultrametric space
$(\mathbb{B}, d)$, where
 $d((x_{n})_{n\in\omega}, (y_{n})_{n\in\omega})= \min\{m: x_{n}=y_{n}, \ n\geq  m\}$. By Proposition 2.1, there are plenty metrically Ramsey ultrafilters on $\mathbb{B}$ in ZFC. Applying Proposition 1.6 (iii), we can find ultrafilters on $\mathbb{B}$ which are not metrically Ramsey.
 \vspace{5 mm}

 3. By [4, Theorem 6.2], there is a coloring $\chi: [\mathbb{R}]^{2}\longrightarrow\{0,1\}$  such that if $X\subset \mathbb{R}$ and
 $[X]^{2}$ is $\chi$-monochrome then $|X|\leq\omega$.

 We endow $\mathbb{R}$ with the natural metric $d(x,y)=|x-y|$ and ask

 \vspace{4 mm}
{\bf Question 4.3.} {\it Does there exist an isometric coloring $\chi: [\mathbb{R}]^{2}\longrightarrow\{0,1\}$  such that if $[X]^{2}$ is monochrome
then $|X|\leq\omega$?}.

We endow the Cantor cube $\{0,1\}^{\omega}$ with the standard metric and ask

\vspace{3 mm}
{\bf Question 4.4.} {\it Does there exist an isometric coloring
$\chi:$ $[\{ 0,1 \}^{\omega}]^{2}\rightarrow\{0,1\}$  such that if $[X]^{2}$ is monochrome
then $|X|\leq\omega$?}.

\bibliography{mybibfile}

\vskip 15pt

CONTACT INFORMATION

I.~Protasov: \\
Faculty of Computer Science and Cybernetics  \\
        Kyiv University  \\
         Academic Glushkov pr. 4d  \\
         03680 Kyiv, Ukraine \\ i.v.protasov@gmail.com

\medskip

K.~Protasova:\\
Faculty of Computer Science and Cybernetics \\
        Kyiv University  \\
         Academic Glushkov pr. 4d  \\
         03680 Kyiv, Ukraine \\ ksuha@freenet.com.ua

\end{document}